\documentclass[11pt, a4paper, reqno]{amsart}
\usepackage[english]{babel}
\usepackage{amsmath}
\usepackage{amsfonts}
\usepackage{amssymb}
\usepackage{url}
\usepackage{amsthm}

\newcommand{\R}{\mathbb{R}}

\newcommand{\N}{\mathbb{N}}

\newcommand{\Ca}{{\rm cap}}
\newcommand{\capa}{{\rm cap}}

\newcommand{\dd}{\, d}
\newcommand{\var}{\varepsilon}
\newcommand{\pois}[1]{}

\DeclareMathOperator*{\diam}{diam}

\DeclareMathOperator*{\dist}{dist}

\def\mvint_#1{\mathchoice%
          {\mathop{\kern 0.2em\vrule width 0.6em height 0.69678ex depth -0.58065ex
                  \kern -0.8em \intop}\nolimits_{\kern -0.4em#1}}%
          {\mathop{\kern 0.1em\vrule width 0.5em height 0.69678ex depth -0.60387ex
                  \kern -0.6em \intop}\nolimits_{#1}}%
          {\mathop{\kern 0.1em\vrule width 0.5em height 0.69678ex depth -0.60387ex
                  \kern -0.6em \intop}\nolimits_{#1}}%
          {\mathop{\kern 0.1em\vrule width 0.5em height 0.69678ex depth -0.60387ex
                  \kern -0.6em \intop}\nolimits_{#1}}}
\def\mvintslides_#1{\mathchoice%
          {\mathop{\kern 0.1em\vrule width 0.5em height 0.697ex depth -0.581ex
                  \kern -0.6em \intop}\nolimits_{\kern -0.4em#1}}%
          {\mathop{\kern 0.1em\vrule width 0.3em height 0.697ex depth -0.604ex
                  \kern -0.4em \intop}\nolimits_{#1}}%
          {\mathop{\kern 0.1em\vrule width 0.3em height 0.697ex depth -0.604ex
                  \kern -0.4em \intop}\nolimits_{#1}}%
          {\mathop{\kern 0.1em\vrule width 0.3em height 0.697ex depth -0.604ex
                  \kern -0.4em \intop}\nolimits_{#1}}}

\def\vint_#1{\mathchoice
          {\mathop{\vrule width 6pt height 3 pt depth -2.5pt
                  \kern -8pt \intop}\nolimits_{#1}}%
          {\mathop{\vrule width 5pt height 3 pt depth -2.6pt
                  \kern -6pt \intop}\nolimits_{#1}}%
          {\mathop{\vrule width 5pt height 3 pt depth -2.6pt
                  \kern -6pt \intop}\nolimits_{#1}}%
          {\mathop{\vrule width 5pt height 3 pt depth -2.6pt
                  \kern -6pt \intop}\nolimits_{#1}}}
\def\sqr#1#2{{\vcenter{\hrule height.#2pt\hbox{\vrule
     width.#2pt height#1pt\kern#1pt\vrule width.#2pt}\hrule height.#2pt}}}

\def\Xint#1{\mathchoice
   {\XXint\displaystyle\textstyle{#1}}%
   {\XXint\textstyle\scriptstyle{#1}}%
   {\XXint\scriptstyle\scriptscriptstyle{#1}}%
   {\XXint\scriptscriptstyle\scriptscriptstyle{#1}}%
   \!\int}
\def\XXint#1#2#3{{\setbox0=\hbox{$#1{#2#3}{\int}$}
     \vcenter{\hbox{$#2#3$}}\kern-.5\wd0}}

\def\vint{\Xint-}

\theoremstyle{plain}
\newtheorem{theorem}[equation]{Theorem}
\newtheorem{lemma}[equation]{Lemma}

\numberwithin{equation}{section}

\theoremstyle{definition}
\newtheorem{definition}[equation]{Definition}
\newtheorem{example}[equation]{Example}

\theoremstyle{remark}
\newtheorem{remark}[equation]{Remark}
\title[Equivalence and of $p$--fatness and Hardy's
inequality]{Equivalence and self--improvement of $p$--fatness and
  Hardy's inequality, and association with uniform perfectness}
\author{Riikka Korte \and Nageswari Shanmugalingam}

\thanks{The first author is supported by the Finnish Academy of
  Science and Letters, Vilho, Yrj\"o and Kalle V\"ais\"al\"a
  Foundation, and the second author is partially supported by the NSF
  grant DMS-0355027. }
\begin{document}



\begin{abstract}
  We present an easy proof that $p$--Hardy's inequality implies
  uniform $p$--fatness of the boundary when $p=n$. The proof works
  also in metric space setting and demonstrates the self--improving
  phenomenon of the $p$--fatness. We also explore the relationship
  between $p$--fatness, $p$--Hardy inequality, and the uniform
  perfectness for all $p\ge 1$, and demonstrate that in the Ahlfors
  $Q$--regular metric measure space setting with $p=Q$, these three
  properties are equivalent. When $p\neq 2$, our results are new even
  in the Euclidean setting.
\end{abstract}

\maketitle

\section{Introduction}

The purpose of this paper is to study the relation between
$p$--Hardy's inequality
\[
\int_\Omega\frac{|u(x)|^p}{\dist(x,\partial\Omega)^p}\dd x\leq
C\int_\Omega |\nabla u(x)|^p\dd x
\]
for all $u\in C_0^\infty(\Omega)$ and $C$ independent of $u$, the
uniform perfectness of $\partial\Omega$, and the uniform $p$--fatness
of $X\setminus\Omega$ in the metric space setting. By $p$--fatness we
mean a capacitary version of the measure thickness condition. Rather
surprisingly, these analytic, metric and geometric conditions turn out
to be equivalent in certain situations. We also consider
self--improving phenomena related to these conditions. Our results are
new even in the Euclidean setting, when $p\neq 2$.

The fact that when $p=n$, a domain satisfies $p$--Hardy's inequality
if and only if the complement is uniformly $p$--fat, was first proved
by Ancona~\cite{ancona86} in $\R^2$. Later, these results were
generalized for all $n=p>1$ by Lewis~\cite{lewis88}. Sugawa proved
in~\cite{sugawa03} that for $n=p=2$ these conditions are equivalent to
the uniform perfectness of the complement in the Euclidean plane.
See also Buckley--Koskela~\cite{buckleykoskela04} for studies relevant
to Orlicz--Sobolev spaces.


In metric spaces, for all $p>1$, it has been shown that uniform
$p$--fatness of the complement of a domain implies that the domain
supports $p$--Hardy's inequality under some conditions,
see~\cite{bjornms01}.  See also~\cite{kkm00} for similar results
involving a measure thickness condition. In~\cite{Lehr07}, the
equivalence of the $p$--fatness and a pointwise Hardy's inequality
has been studied. In this paper, we prove that if a metric space is
Ahlfors $Q$--regular and satisfies a weak $(1,Q)$--Poincar\'e
inequality, then the support of a $Q$--Hardy inequality on a domain
implies uniform $Q$--fatness of the complement of the domain.  Our
proof is rather transparent and it is based on estimating the
Hausdorff--content of the boundary.

We will also prove a self--improvement property for both uniform
$Q$--fatness and $Q$--Hardy's inequality in the setting of Ahlfors
$Q$--regular metric measure spaces. That is, if a set satisfies
$Q$--Hardy's inequality or is uniformly $Q$--fat, then there exists
$q<Q$ such that the set satisfies $q$--Hardy's inequality or is
uniformly $q$--fat, respectively.  The self--improving property of
Hardy's inequality has been studied in~\cite{koskelazhong03} and that
of uniform $p$--fatness in~\cite{bjornms01} and in~\cite{mikkonen96}.
Our approach gives a more elementary proof of self--improvement of
uniform $p$--fatness when $p=Q$.

\section{Preliminaries}

We assume that $X=(X,d,\mu)$ is a metric measure space equipped with a
metric $d$ and a Borel regular outer measure $\mu$ such that
$0<\mu(B)<\infty$ for all balls $B=B(x,r)=\{y\in X\,:\, d(x,y)<r\}$.
The measure $\mu$ is said to be \emph{doubling} if there exists a
constant $c_D\geq 1$, called the \emph{doubling constant}, such that
\[
\mu(B(x,2r))\leq c_D\mu(B(x,r))
\]
for all $x\in X$ and $r>0$. The measure is \emph{$Q$--regular} if
there exists a constant $c_A\geq 1$ such that
\[
\frac1{c_A}r^Q\leq\mu(B(x,r))\leq c_A r^Q
\]
for all $x\in X$ and $0<r<\diam(X)$. The $n$--dimensional Lebesgue
measure on $\R^n$ is $n$--regular. The Hausdorff $s$--content of
$E\subset X$ is
\begin{equation}\label{eqn:hausdorff_content}
  \mathcal H^s_\infty(E)=\inf\sum_{i\in I}r_i^s,
\end{equation}
where the infimum is taken over all countable covers
$\{B(x_i,r_i)\}_{i\in I}$ of $E$, with each $B(x_i,r_i)\cap E$
non-empty. In addition, we may assume that $x_i\in E$ for every $i\in
I$, because that may increase the Hausdorff content at most by a
multiplicative factor $2^s$.

A non-negative Borel measurable function $g_u$ on $X$ is said to be a
\emph{$p$--weak upper gradient} of a function $u$ on $X$ if there is a
non-negative Borel measurable function $\rho\in L^p(X)$ such that for
all rectifiable curves $\gamma$ in $X$, denoting the end points of
$\gamma$ by $x$ and $y$, we have either
\[
  |u(x)-u(y)|\le \int_\gamma g_u\, ds,
\]
or $\int_\gamma \rho\, ds=\infty$.  Let $1\leq p<\infty$. If $u$ is a
function that is integrable to power $p$ in $X$, let
\[
\|u\|_{N^{1,p}(X)}=\left(\int_X |u|^p\dd\mu+\inf_{g_u}\int_X
  g_u^p\dd\mu\right)^\frac1p,
\]
where the infimum is taken over all $p$--weak upper gradients of $u$.
The \emph{Newtonian space} on $X$ is the quotient space
\[
N^{1,p}(X)=\{u\,:\,\|u\|_{N^{1,p}(X)}<\infty\}/\sim,
\]
where $u\sim v$ if and only if $\|u-v\|_{N^{1,p}(X)}=0$,
see~\cite{shanmugalingam00}. We
define $N^{1,p}_0(\Omega)$ to be the set of functions $u\in
N^{1,p}(\Omega)$ that can be extended to $N^{1,p}(X)$ so that the
extensions are zero on $X\setminus\Omega$ $p$--quasieverywhere.

We say that $X$ supports a \emph{weak $(1,p)$--Poincar\'e inequality}
if there exist constants $c_p>0$ and $\tau\geq 1$ such that for all
balls $B(x,r)$ of $X$, all locally integrable functions $u$ on $X$ and
for all $p$--weak upper gradients $g_u$ of $u$, we have
\begin{equation}\label{eqn:poincare}
  \mvint_{B(x,r)} |u-u_{B(x,r)}|\dd\mu\leq c_p r\left( \mvint_{B(x,\tau r)}g_u^p\dd\mu\right)^\frac{1}{p},
\end{equation}
where
\begin{equation*}
  u_B=\mvint_Bu\dd\mu=\frac{1}{\mu(B)}\int_B u\dd\mu.
\end{equation*}

We point out here that if $X$ is the Euclidean space $\R^n$ equipped
with the $n$--dimensional Lebesgue measure and the Euclidean metric,
then $N^{1,p}(X)=W^{1,p}(\R^n)$, the classical Sobolev space.
Moreover, $\R^n$ supports a weak $(1,1)$--Poincar\'e inequality.
\begin{definition} Let $\Omega$ be an open set in $X$ and $E$ be a
  closed subset of $\Omega$. The \emph{$p$--capacity} of $E$ with
  respect to $\Omega$ is
  \[
  \Ca_p(E,\Omega)=\inf\int_X g_u^p\dd\mu,
  \]
  where the infimum is taken over all functions $u$ with $p$--weak
  upper gradients $g_u$ such that $u_{|E}=1$ and $u_{|X\setminus
    \Omega}=0$.  Should there be no such function $u$, then
  $\Ca_p(E,\Omega)=\infty$.
\end{definition}

A metric space $X$ is said to be \emph{linearly locally connected}
(LLC) if there is a constant $C\geq 1$ so that for each $x\in X$ and
$r>0$, the following two conditions hold:
\begin{enumerate}
\item any pair of points in $B(x,r)$ can be joined in $B(x,Cr)$,
\item any pair of points in $X\setminus\overline B(x,r)$ can be joined
  in $X\setminus \overline B(x,r/C)$.
\end{enumerate}
By \emph{joining} we mean joining by a path.  Note that if a complete
$Q$--regular space, $Q>1$, supports a weak $(1,Q)$--Poincar\'e
inequality, then it satisfies the LLC--condition, see for
example~\cite{heinonenkoskela98} or~\cite{korte06}.

\begin{definition}\label{def:up}
  We say that a set $E\subset X$ is \textit{uniformly perfect} if $E$
  is not a singleton set, and there is a constant $c_{UP}\geq 1$ so
  that for each $x\in E$ and $r>0$ the set $E\cap
  B(x,c_{UP}r)\setminus B(x,r)$ is nonempty whenever the set
  $E\setminus B(x,c_{UP}r)$ is nonempty.
\end{definition}

For more information about uniform perfectness, see for
example~\cite{jarvivuorinen96} and~\cite{sugawa03}.

A set $E\subset X$ is said to be \textit{uniformly $p$--fat} 
if there exists a constant $c_0>0$ so that for 
every point $x\in E$ and for all $0<r<\infty$,
\begin{equation}\label{eqn:p-fat}
  \frac{\Ca_p(B(x,r)\cap E,B(x,2r))}{\Ca_p(B(x,r),B(x,2r))}
  \geq c_0.
\end{equation}
This condition is stronger than the Wiener criterion.  Uniform
$p$--fatness is a capacitary version of the uniform measure thickness
condition, see for example~\cite{kkm00}.

\begin{definition}\label{def:hardy}
  Let $1<p<\infty$. The set $\Omega\subset X$ satisfies $p$--Hardy's
  inequality if there exists $0<c_H<\infty$ such that for all $u\in
  N^{1,p}_0(\Omega)$,
  \begin{equation}\label{eqn:hardy}
    \int_\Omega \left(\frac{|u(x)|}{\dist(x,X\setminus\Omega)}\right)^p\dd\mu(x)\leq c_H\int_\Omega g_u(x)^p\dd\mu(x).
  \end{equation}
  Here $g_u$ is a $p$--weak upper gradient of $u$. Here we use
  $\dist(x,X\setminus\Omega)$ instead of $\dist(x,\partial\Omega)$
  since in general the latter quantity can be larger than the former
  one.
\end{definition}
Hardy's inequality has been studied for example
in~\cite{davies95},~\cite{hajlasz99},~\cite{lewis88},
and~\cite{tidblom04}.  Hardy's inequality has been used also to
characterize Sobolev functions with zero boundary values,
see~\cite{kinnunenmartio97} and~\cite{lewis88}.
\section{Main results}

In this section, we show that $Q$--Hardy's inequality on $\Omega$
implies uniform $Q$--fatness of the complement. Our method also shows
that $Q$--fatness is a self--improving property. To simplify notation,
we will assume $X$ to be unbounded throughout this section. However,
for our arguments, it is immaterial what happens outside $\overline
\Omega$, and therefore our arguments work also if $X$ is bounded,
provided we adjust the conditions of uniform perfectness and uniform
fatness to the bounded setting. Notice also that if $\Omega$ is a
domain, then $X\setminus\Omega$ can be replaced by $\partial \Omega$
in our arguments.

\begin{theorem}\label{theorem:hardy_fat}Let $(X,d,\mu)$ be a complete
  $Q$--regular metric measure space supporting a weak
  $(1,Q)$--Poincar\'e inequality, and $\Omega\subset X$ be an open
  subset. If $\Omega$ satisfies $Q$--Hardy's inequality, then
  $X\setminus \Omega$ is uniformly $(Q-\varepsilon)$--fat for some
  $\varepsilon>0$.
\end{theorem}

We split the proof into two parts. First in
Lemma~\ref{lemma:hardy_up}, we show that Hardy's inequality implies
uniform perfectness of the complement. Then in
Theorem~\ref{theorem:up_fat}, we show that uniform perfectness implies
$(Q-\varepsilon)$--fatness with some $\varepsilon>0$. Recall that we
assume $X$ to be unbounded.

\begin{lemma}\label{lemma:hardy_up}
  Let $X$ be as in Theorem~\ref{theorem:hardy_fat}. If $\Omega\subset
  X$ satisfies $Q$--Hardy's inequality, then $X\setminus \Omega$ is
  uniformly perfect and unbounded.
\end{lemma}
\begin{proof}
  Fix $m>4$ and suppose that $\Omega$ satisfies Hardy's
  inequality~\eqref{eqn:hardy} and that $X\setminus\Omega$ is not
  uniformly perfect with respect to the constant $m$ or that
  $X\setminus\Omega$ is bounded. In both cases, there exists $x_0\in
  X\setminus\Omega$ and $r_0>0$ such that $B(x_0,mr_0)\setminus
  B(x_0,r_0)\subset \Omega$. We will deduce an upper bound for such
  $m$ independent of $x_0$ and $r_0$, and hence conclude that
  $X\setminus\Omega$ is uniformly perfect for any constant larger than
  this upper bound and that $X\setminus\Omega$ cannot be bounded.

  Define $u:X\rightarrow [0,\infty)$ so that
  \begin{equation*}
    u(x)=
    \begin{cases}
      \left(\frac{d(x_0,x)}{r_0}-1\right)_+, & d(x_0,x)\leq 2r_0,\\
      1, & 2r_0<d(x_0,x)<\frac{m\,r_0}{2},\\
      \left(2-\frac{2d(x_0,x)}{m\,r_0}\right)_+, &
      \frac{m\,r_0}{2}\leq d(x_0,x).
    \end{cases}
  \end{equation*}
  Now the minimal upper gradient of $u$ satisfies
  \begin{equation}\label{eqn:h_o}
    \int_\Omega g_u^Q\dd\mu 
    \leq \left(\frac{1}{r_0}\right)^Q\mu(B(x_0,2r_0))+\left(\frac{2}{m r_0}\right)^Q\mu(B(x_0,m r_0))
    \leq c_A 2^{Q+1}.
  \end{equation}
  Next, we show that
  \begin{equation}\label{eqn:h_v}
    \int_\Omega \frac{u(x)^Q}{\dist(x,X\setminus \Omega)^Q}\dd\mu(x)\geq c \log(m/4),
  \end{equation}
  where $c>0$ is a constant that depends only on $c_A$ and $Q$.  For
  $x\in X$ and $0<r<R$, we denote the annulus
  $A(x,r,R)=B(x,R)\setminus B(x,r)$.  Let $n\in\N$ be the unique
  number such that $2^n\le m<2^{n+1}$. Since $m>4$, we have $n\ge 2$.
  Then
  \[
  A(x_0,2r_0,mr_0/2)\supset \bigcup_{k=1}^{n-1}A(x_0,2^kr_0,2^{k+1}r_0).
  \]
  As $X$ is quasiconvex (which follows from the Poincar\'e inequality,
  see for example~\cite{korte06}) and hence path-connected, and as
  $X\setminus B(x_0,2^{k+1}r_0)$ is non-empty, there is a point
  $y_k\in A(x_0,2^kr_0,2^{k+1}r_0)$ such that $d(x_0,y_k)=\dfrac32
  2^k$; hence the ball $B(y_k,2^{k-1}r_0)\subset
  A(x_0,2^kr_0,2^{k+1}r_0)$. Thus
  \begin{align*}
    \int_\Omega\frac{u(x)^Q}{\text{dist}(x,X\setminus\Omega)^Q}\, d\mu
    & \ge
    \int_{A(x_0,2r_0,mr_0)}\frac{1}{d(x_0,x)^Q}\, d\mu \\
    &\ge \sum_{k=1}^{n-1} \int_{A(x_0,2^kr_0,2^{k+1}r_0)}\frac{1}{d(x_0,x)^Q}\, d\mu\\
    &\ge \sum_{k=1}^{n-1}\int_{B(y_k,2^{k-1}r_0)}\frac{1}{d(x_0,x)^Q}\, d\mu\\
    &\ge \sum_{k=1}^{n-1}\frac{1}{(2^{k+1}r_0)^Q}\, \mu(B(y_k,2^{k-1}r_0))\\
    &\ge \frac{n-1}{4^Q\, c_A}.
  \end{align*}
  Since $n> \dfrac{\log(m/2)}{\log(2)}$, we see that $n-1 > \dfrac{\log(m/4)}{\log(2)}$. Thus,
  \[
  \int_\Omega\frac{u(x)^Q}{\text{dist}(x,X\setminus\Omega)^Q}\, d\mu
  > \frac{\log(m/4)}{4^Q\, c_A\, \log(2)}=c\, \log(m/4).
  \]
  By combining~\eqref{eqn:h_o} and~\eqref{eqn:h_v}, and the fact that
  $u$ satisfies Hardy's inequality~\eqref{eqn:hardy}, it follows
  that
  \begin{equation*}
    c\log(m/4)< 2^{Q+1}c_H\, c_A.
  \end{equation*}
  Hence $m< 4\exp(2^{Q+1}c_H\, c_A/c)$, and
  therefore $X\setminus \Omega$ is uniformly perfect with constant
  $c_{UP}= 4\exp(2^{Q+1}c_Hc_A/c)$ and $X\setminus\Omega$ is
  unbounded.
\end{proof}

The following example shows that $p$--Hardy's inequality with $p\neq
Q$ does not imply uniform perfectness.

\begin{example}\label{example:1}
  If $X=\R^n$, $1<p<n$, and $\Omega=B(0,1)\setminus\{0\}$, then
  $\Omega$ supports $p$--Hardy's inequality even though
  $X\setminus\Omega$ is neither uniformly perfect nor uniformly
  $p$--fat, see~\cite[p. 179]{lewis88}. When $p>n$, even single points
  have positive $p$--capacity and hence $X\setminus \Omega$ is
  uniformly $p$--fat and supports $p$--Hardy's inequality but
  $X\setminus\Omega$ is not uniformly perfect.
\end{example}

The uniform perfectness of the boundary implies uniform $q$--fatness
of the complement for all $q>Q-\varepsilon$. We get a quantitative
estimate for $\varepsilon>0$ that depends only on $c_{UP}$.

\begin{theorem}\label{theorem:up_fat}
  Let $(X,d,\mu)$ be a complete $Q$--regular metric measure space.
  Suppose that $X$ supports a weak $(1,Q)$--Poincar\'e inequality.
  Let $\Omega\subset X$ be an open subset.  If $X\setminus \Omega$ is
  uniformly perfect and unbounded, then there exists a constant
  $\varepsilon>0$ such that $X\setminus \Omega$ is uniformly
  $(Q-\varepsilon)$--fat.
\end{theorem}
We begin the proof with an elementary inequality.
\begin{lemma}\label{lemma:epsilon_C}
  For every $C>0$ there exists $0<\varepsilon_C<1$ such that for all
  $0<\varepsilon<\varepsilon_C$ and $a,b>0$,
  \[
  a^\varepsilon+b^\varepsilon\geq\left(a+b+C\min\{a,b\}\right)^\varepsilon.
  \]
\end{lemma}

\begin{proof}
  We may assume that $a\geq b=1$. Therefore, it is enough to prove
  that
  \[
  (a+C+1)^\varepsilon-a^\varepsilon\leq 1
  \]
  when $0<\varepsilon <1$ is sufficiently small and $a\geq1$. As
  $f(a)=(a+C+1)^\varepsilon-a^\varepsilon$ is a decreasing function,
  and hence for every $a\geq 1$
  \[
  (a+C+1)^\varepsilon-a^\varepsilon\leq (1+C+1)^\varepsilon-1^\varepsilon,
  \]
  it is enough to choose $\varepsilon$ so that
  \[
  \varepsilon\leq \varepsilon_C= \frac{\log 2}{\log(C+2)}.
  \]
\end{proof}
In the proof of Theorem~\ref{theorem:up_fat}, we first obtain an
estimate for the Hausdorff--content of the boundary. Then the
following result is needed to get capacitary estimates. For a proof,
see Theorem 5.9 in~\cite{heinonenkoskela98}.

\begin{lemma}\label{seuraus}
  Suppose that $(X,d,\mu)$ is a $Q$--regular space. Suppose further
  that $X$ admits a weak $(1,p)$--Poincar\'e inequality for some
  $1\leq p\leq Q$.  Let $E\subset B(x,r)$ be a compact set.  If
  \[
  \mathcal H_\infty^s(E)\geq \lambda r^s
  \]
  for some $s> Q-p$ and $\lambda>0$, then
  \[
  \Ca_p(E,B(x,2r))\geq c\lambda \Ca_p(B(x,r),B(x,2r)).
  \]
  The constant $c$ depends only on $s$ and on the data associated with
  $X$.
\end{lemma}

\begin{proof}[Proof of Theorem~\ref{theorem:up_fat}]
  Let $\Omega\subset X$ be open, $X\setminus \Omega$ uniformly perfect
  with constant $c_{UP}>1$, and $\alpha>1$. Fix $x_0\in X\setminus
  \Omega$ and $r_0>0$. Let $A=\overline B(x_0,r_0)\setminus \Omega$,
  and $0<\varepsilon<\varepsilon_{\alpha c_{UP}}$, where
  $\varepsilon_{\alpha c_{UP}}$ is as in Lemma~\ref{lemma:epsilon_C}.
  First we estimate the Hausdorff $\varepsilon$--content of $A$.  Let
  $\mathcal F$ be a family of balls covering $A$.  Because $A$ is
  compact, we may assume that $\mathcal F$ consists of a finite number
  of balls. We may also assume that all the balls in $\mathcal F$ are
  centered at $A$.

  If there exists balls $B(x_i,r_i)$ and $B(x_j,r_j)$ in $\mathcal F$
  such that
  \begin{equation}\label{eqn:r_alpha}
    r_i\leq \alpha r_j
  \end{equation} 
  and 
  \begin{equation}\label{eqn:pallot_lahella}
    B(x_i,c_{UP}r_i)\cap B(x_j,r_j)\neq\emptyset,
  \end{equation}
  then (when $r_j\leq r_i$)
  \[
  \left(B(x_i,r_i)\cup B(x_j,r_j)\right)\subset 
  B\left(x_i,r_i+r_j+\alpha c_{UP}\min\{r_i,r_j\}\right)
  \]
  or (when $r_i\leq r_j$)
  \[
  \left(B(x_i,r_i)\cup B(x_j,r_j)\right)\subset 
  B\left(x_j,r_i+r_j+\alpha c_{UP}\min\{r_i,r_j\}\right),
  \]
  and by Lemma~\ref{lemma:epsilon_C},
  \[
  r_i^\varepsilon+r_j^\varepsilon\geq 
  (r_i+r_j+\alpha c_{UP}\min\{r_i,r_j\})^\varepsilon.
  \] 
  Thus, we may replace balls 
  $B(x_i,r_i)$ and $B(x_j,r_j)$ with 
  \[
  B(x_i , r_i+ r_j + \alpha c_{UP} \min\{r_i,r_j\})\textrm{ or } 
  B(x_j , r_i+ r_j + \alpha c_{UP} \min\{r_i,r_j\})
  \]
  in the covering $\mathcal F$ so that the sum
  \[
  \sum_{B(x,r)\in \mathcal F}r^\varepsilon
  \]
  does not increase. We continue this process until there is no pair
  of balls satisfying~\eqref{eqn:pallot_lahella}
  and~\eqref{eqn:r_alpha}.  Because the number of balls in $\mathcal
  F$ decreases in each step and $\mathcal F$ consists of a finite
  number of balls, the process ends after a finite number of
  replacements.

  Let $B(x_1,r_1)\in\mathcal F$ be the ball containing $x_0$. Because
  $X\setminus \Omega$ is uniformly perfect and
  unbounded, 
  the set
  \[
  \left(B(x_1,c_{UP}\,r_1)\setminus B(x_1,r_1)\right)\cap (X\setminus\Omega)
  \]
  is nonempty. Now there are two possibilities: either 
  \[
  (B(x_1,c_{UP}\,r_1)\setminus B(x_1,r_1))\cap 
  \left(X\setminus \overline B(x_0,r_0)\right)\neq \emptyset
  \]
  or 
  \[
  (B(x_1,c_{UP}\,r_1)\setminus B(x_1,r_1))\cap A\neq \emptyset,
  \]
  because $(X\setminus \Omega)\subset A \cup (X\setminus \overline
  B(x_0,r_0))$.  In the latter case, there exists $B(x_2,r_2)\in
  \mathcal F$ such that $B(x_1,r_1)\neq B(x_2,r_2)$ and
  \[
  B(x_2,r_2)\cap B(x_1,c_{UP}\,r_1)\neq\emptyset,
  \]
  because $\mathcal F$ covers $A$. Now the balls $B(x_1,r_1)$ and
  $B(x_2,r_2)$ satisfy condition~\eqref{eqn:pallot_lahella}.
  Hence~\eqref{eqn:r_alpha} fails, that is, $r_2<r_1/\alpha$.

  We continue inductively in the same way: For a ball $B(x_i,r_i)\in
  \mathcal F$, either
  \[
  B(x_i,c_{UP}r_i)\cap (X\setminus \overline B(x_0,r_0))\neq \emptyset
  \]
  or there exists a ball $B(x_{i+1},r_{i+1})\in\mathcal F$ such that
  $r_{i+1}\leq r_i/\alpha$ and $B(x_i,r_i)$ and $B(x_{i+1},r_{i+1})$
  satisfy the condition~\eqref{eqn:pallot_lahella}.

  Thus we obtain a chain of distinct balls
  $\{B(x_i,r_i)\}_{i=1}^n\subset \mathcal F$ such that
  $r_i\leq\alpha^{1-i}r_1$, (since $r_i\leq r_{i-1}/\alpha$, we have
  $B(x_i,r_i)\neq B(x_j,r_j)$ if $i\neq j$).
  \[ 
  B(x_i,c_{UP}\,r_i)\cap B(x_{i+1},r_{i+1})\neq\emptyset 
  \]
  for every $i=1,\ldots,n-1$, and 
  \[
  B(x_n,c_{UP}\,r_n)\cap \left(X\setminus \overline 
    B(x_0,r_0)\right)\neq \emptyset.
  \]
  It follows that
  \[
  r_0\leq \sum_{i=1}^n(c_{UP}+1)r_i\leq
  (c_{UP}+1)\frac{\alpha}{\alpha-1}\,r_1
  \]
  and we have a lower bound for $r_1$:
  \[
  r_1\geq \frac{\alpha-1}{\alpha(c_{UP}+1)}\, r_0.
  \]

  We may choose $\alpha=2$ and thus
  \[
  \sum_{B(x,r)\in\mathcal F}r^\varepsilon\geq r_1^\varepsilon\geq
  \frac{1}{(2c_{UP}+2)^\varepsilon}r_0^\varepsilon.
  \]
  By~\cite{keithzhong06}, there exists $\varepsilon>0$ such that $X$
  satisfies a weak $(1,Q-\varepsilon)$--Poincar\'e inequality. Fix
  such an $\varepsilon<\varepsilon_{\alpha c_{UP}}$.  Now by
  Lemma~\ref{seuraus},
  \[
  \Ca_p(B(x_0,r_0)\setminus \Omega,B(x_0,2r_0))\geq
  c\,\Ca_p(B(x_0,r_0),B(x_0,2r_0))
  \]
  for every $Q-\varepsilon<p\leq Q$, where $c$ depends on
  $\varepsilon$ and $c_{UP}$, but is independent of $x_0$ and $r_0$.
\end{proof}
It is known that uniform $p$--fatness is a self--improving phenomenon,
see~\cite{bjornms01}.
\begin{theorem}
  Let $X$ be a proper linearly locally convex metric space endowed
  with a doubling Borel regular measure supporting a weak
  $(1,q_0)$--Poincar\'e inequality for some $q_0$ with $1\leq
  q_0<\infty$. Let $p>q_0$ and suppose that $E\subset X$ is uniformly
  $p$--fat. Then there exists $q<p$ so that $E$ is uniformly $q$--fat.
\end{theorem}

\begin{remark}\label{remark:easier}
  The proof of Theorem~\ref{theorem:up_fat} gives a new and easier
  proof for the self--improvement when $p=Q$.
\end{remark}
\begin{remark}
  To complete the picture, note that uniform $p$--fatness for any
  $p\leq Q$ implies uniform perfectness.  To see this, suppose that
  $X$ supports a $(1,p)$--Poincar\'e inequality for some $1\le p<Q$,
  and that $X\setminus\Omega$ is uniformly $p$--fat. We will show that
  $X\setminus\Omega$ is uniformly perfect. Fix $x_0\in\partial\Omega$
  and $0<r<1$. Suppose that $B(x_0,r)\setminus
  B(x_0,r/m)\subset\Omega$ for some $m>1$. We will demonstrate that
  $m$ has an upper bound that is independent of $x_0$ and $r$. Indeed,
  \[
  \Ca_p(B(x_0,r), B(x_0,2r))\ge \frac1C\, r^{Q-p}.
  \]
  Also, as the function 
  \[
  g(x)= \frac{1}{\log(r/\rho)}\,
  \frac{1}{d(x_0,x)}\chi_{B(x_0,r)\setminus B(x_0,\rho)}
  \]
  is an upper gradient of the function
  \[
  u(x)=\min\left\{
    1,\max\left\{0,\frac{\log(d(x_0,x)/\rho)}{\log(r/\rho)}\right\}\right\},
  \]
  with $u=0$ on $B(x_0,\rho)$ and $u=1$ on $X\setminus B(x_0,r)$; hence
  \begin{align*}
    \Ca_p(B(x_0,r)\setminus\Omega,B(x_0,2r))&
    \le \Ca_p(B(x_0,r/m),B(x_0,2r))\\
    &\le \frac{C}{\log(m)^p}\, r^{Q-p}.
  \end{align*}
The last estimate can be proved in the same way as in the proof of
  Lemma~\ref{seuraus}. We have by uniform $p$--fatness of
  $X\setminus\Omega$ that
  \begin{align*}
    \frac{C}{\log(m)^p}\, r^{Q-p}&\ge \Ca_p(B(x_0,r)\setminus\Omega,B(x_0,2r))\\
    &\ge \frac{1}{c_0} \Ca_p(B(x_0,r), B(x_0,2r)) \ge \frac1{c_0C}\, r^{Q-p},
  \end{align*}
  where $c_0$ is the uniform fatness constant, and therefore $m\le
  e^{C}$.  Thus $X\setminus \Omega$ is uniformly perfect.
\end{remark}

\begin{remark}\label{remark:sharp}
  In the proof of Theorem~\ref{theorem:up_fat}, we need to assume that
  the space supports a weak $(1,Q-\varepsilon)$--Poincar\'e
  inequality. This follows by~\cite{keithzhong06} for some positive
  $\varepsilon$ if the space supports a $(1,Q)$--Poincar\'e
  inequality. However, if we assume a priori the stronger Poincar\'e
  inequality, then our proof gives a quantitative estimate for
  $\varepsilon$. More precisely, if
  \[
  \max\left\{Q-\dfrac{\log(2)}{\log(3)},1\right\}<p<Q
  \]
  and $X$ supports a $(1,p)$--Poincar\'e inequality, then there exists
  $c_p>1$ such that whenever $X\setminus \Omega$ is uniformly perfect
  for some uniform perfectness constant $1\leq c_{UP}<c_p$, then
  $X\setminus \Omega$ is uniformly $p$--fat and hence $\Omega$
  supports a $p$--Hardy inequality by Theorem~\ref{theorem:fat_hardy}.
  The proof of Theorem~\ref{theorem:up_fat} implies the claim if
  \[
  p>Q-\frac{\log(2)}{\log(\alpha c_{UP}+2)}
  \]
  with some $\alpha>1$. So it is enough to have $c_{UP}< c_p=
  2^\frac{1}{Q-p}-2$.  By the assumption on $p$, it is clear that
  $c_p>1$.
\end{remark}
The following examples illustrate the sharpness of
Remark~\ref{remark:sharp}.
\begin{example} If $1<p<Q$, there is a Cantor set $E_p\subset \R^n$
  such that $\Ca_p(E_p)=0$, see~\cite[p. 40]{heinonenkoskela98}. Thus
  the domain $\R^n\setminus E_p$ has uniformly perfect complement,
  which is not uniformly $p$--fat.
\end{example}
\begin{example}If $1\leq p<Q-1$, then any rectifiable curve $\gamma$
  in $X$ is of zero $p$--capacity. In this case, with
  $\Omega=X\setminus\gamma$, we have that $X\setminus\Omega$ is
  uniformly perfect with constant $c_{UP}=1$, but it is not uniformly
  $p$--fat.
\end{example}

The following theorem shows that Hardy's inequality follows from
uniform fatness for all $1<p\leq Q$, see Corollary 6.1
in~\cite{bjornms01}. Note that the LLC--condition is not a serious
restriction in our case, since it follows from the $(1,Q)$--Poincar\'e
inequality, see for example~\cite{korte06}.

\begin{theorem}\label{theorem:fat_hardy}
  Let $X$ be a proper LLC metric space endowed with a doubling Borel
  regular measure supporting a weak $(1,p)$--Poincar\'e inequality,
  and suppose that $\Omega$ is a bounded open set in $X$ with
  $X\setminus\Omega$ uniformly $p$--fat. Then $\Omega$ satisfies
  $p$--Hardy's inequality.
\end{theorem}

The converse of Theorem~\ref{theorem:fat_hardy} is not true in
general, see Example~\ref{example:1}.  As a corollary of
Lemma~\ref{lemma:hardy_up} and Theorems~\ref{theorem:up_fat}
and~\ref{theorem:fat_hardy}, we obtain the following result.  Note
that uniform $p$--fatness implies uniform $q$--fatness for all $q>p$.
\begin{theorem}\label{corollary:main}
  Let $(X,d,\mu)$ be a complete $Q$--regular metric measure space with $Q>1$.
  Suppose that $X$ supports a weak $(1,Q)$--Poincar\'e inequality.
  Let $\Omega\subset X$ be a bounded open subset. Then the following
  conditions are quantitatively equivalent.
  \begin{enumerate}
  \item[(1)] $\Omega$ satisfies $Q$--Hardy's inequality.
  \item[(2)] $X\setminus \Omega$ is uniformly perfect.
  \item[(3)] $X\setminus\Omega$ is uniformly $Q$--fat.
  \item[(4)] $X\setminus\Omega$ is uniformly $(Q-\varepsilon)$--fat
    for some $\var>0$.
  \end{enumerate}
\end{theorem}

Theorem~\ref{theorem:fat_hardy} is stated only for bounded sets but
the proof works also in the unbounded setting. Hence
Theorem~\ref{corollary:main} holds also when $\Omega$ is unbounded if
we require additionally that $X\setminus\Omega$ is unbounded in
conditions (2) and (3).
\section{Maz$'$ya type characterization}

In this section, we present one more characterization of open sets
that is equivalent with the Hardy's inequality. For more information
about this kind of characterizations, see Chapter 2.3
in~\cite{mazja85}.

\begin{theorem}
  Let $X$ be a complete metric space endowed with a doubling measure
  and supporting a weak $(1,p)$--Poincar\'e inequality.  Let $1<p\leq
  Q$. Then $\Omega\subset X$ satisfies $p$--Hardy's inequality if and
  only if for every $K\subset\subset\Omega$, we have
  \begin{equation}\label{eqn:mazja}
    \int_K \dist(x,X\setminus\Omega)^{-p}\dd\mu(x)\leq c\,\Ca_p(K,\Omega).
  \end{equation}
\end{theorem}

\begin{proof}
  First assume that $\Omega$ satisfies $p$--Hardy's inequality.  Let
  $u\in N^{1,p}_0(\Omega)$ such that $u=1$ in $K$.  Then
  \[
  \int_K\dist(x,X\setminus\Omega)^{-p}\dd\mu(x)\leq \int_\Omega
  \frac{|u(x)|^p}{\dist(x,X\setminus\Omega)^p}\dd\mu(x)\leq
  c_H\int_\Omega g_u^p\dd\mu.
  \]

  By taking infimum over all such functions $u$, we
  obtain~\eqref{eqn:mazja}.

  Now assume that equation~\eqref{eqn:mazja} is satisfied. We will
  first prove the claim for Lipschitz--functions that
  have compact support in $\Omega$. By Theorems 2.12 and 4.8
  in~\cite{shanmugalingam01}, such functions form a dense subclass of
  $N_0^{1,p}(\Omega)$, and thus we get the result for all functions in
  $N_0^{1,p}(\Omega)$.

  Let $u\in N_0^{1,p}(\Omega)$ be compactly supported Lipschitz 
function, and denote
  \[
  E_k=\{x\in \Omega\,:\, |u(x)|>2^k\},\, k=1,2,\ldots.
  \]
  Thus by~\eqref{eqn:mazja}, we have
  \begin{align*}
    \int_\Omega\frac{|u(x)|^p}{\dist(x,X\setminus\Omega)^p}\dd\mu(x)
    \leq & \sum_{k=-\infty}^\infty 2^{(k+1)p}\int_{E_k\setminus E_{k+1}}\frac1{\dist(x,X\setminus\Omega)^p}\dd\mu(x)\\
    \leq &\, c\sum_{k=-\infty}^\infty 2^{(k+1)p}\Ca_p(\overline E_{k+1},\Omega)\\
    \leq &\, c\sum_{k=-\infty}^\infty 2^{(k+1)p}\Ca_p(\overline E_{k+1},E_{k}).
  \end{align*}
  Let
  \begin{equation*}
    u_k=
    \begin{cases}
      1, & \textrm{when }|u|\geq 2^{k+1},\\
      \frac{|u|}{2^{k}}-1, & \textrm{when } 2^{k}<|u|<2^{k+1},\\
      0,&\textrm{when }|u|\leq 2^{k}. 
    \end{cases}
  \end{equation*}
  Then $u_k=1$ in $\overline E_{k+1}$ and $u_k=0$ in $X\setminus E_{k}$.
  Therefore,
  \[
  \capa_p(\overline E_{k+1},E_{k})\leq \int_{E_{k}\setminus E_{k+1}} g_{u_k}^p\dd\mu
  \leq 2^{-pk}
  \int_{E_{k}\setminus E_{k+1}} g_u^p\dd\mu.
  \]
  Consequently,
  \begin{align*}
    c\sum_{k=-\infty}^\infty 2^{(k+1)p}\Ca_p(\overline E_{k+1},E_{k})
    \leq &\, c2^{p}\sum_{k=-\infty}^\infty\int_{E_{k}\setminus E_{k+1}}g_u^p\dd\mu\\
    =   &\, c2^{p}\int_\Omega g_u^p\dd\mu,
  \end{align*}
  and the claim follows with $c_H=2^{p}c$.
\end{proof}
\subsubsection*{Acknowledgements}
The authors would like to thank Juha Kinnunen for introducing the
subject and for valuable comments.
\def\cprime{$'$} \def\cprime{$'$} \def\cprime{$'$} \def\cprime{$'$}
  \def\cprime{$'$} \def\ocirc#1{\ifmmode\setbox0=\hbox{$#1$}\dimen0=\ht0
  \advance\dimen0 by1pt\rlap{\hbox to\wd0{\hss\raise\dimen0
  \hbox{\hskip.2em$\scriptscriptstyle\circ$}\hss}}#1\else {\accent"17 #1}\fi}
  \def\cprime{$'$}

\vspace{0.5cm}
\noindent
\small{\textsc{R.K.}}\\
\small{\textsc{Institute of Mathematics},}\\
\small{\textsc{P.O. Box 1100},}\\
\small{\textsc{FI-02015 Helsinki University of Technology},}\\
\small{\textsc{Finland}}\\
\footnotesize{\texttt{rkorte@math.hut.fi}}

\vspace{0.3cm}
\noindent
\small{\textsc{N.S.}}\\
\small{\textsc{Department of Mathematical Sciences},}\\
\small{\textsc{P.O. Box 210025},}\\
\small{\textsc{University of Cincinnati},}\\
\small{\textsc{Cincinnati, OH 45221-0025}}\\
\small{\textsc{U.S.A.}}\\
\footnotesize{\texttt{nages@math.uc.edu}}


\end{document}